\documentclass[11pt]{article}

\usepackage[letterpaper,margin=1in]{geometry}
\usepackage{amsmath, amssymb, amsthm}
\usepackage[T1]{fontenc}
\usepackage{lmodern}
\usepackage{microtype}
\usepackage{graphicx}
\usepackage{booktabs,array}
\usepackage{xspace}
\usepackage[dvipsnames]{xcolor}
\usepackage[numbers,sort&compress]{natbib}
\usepackage{xurl}
\usepackage[colorlinks=true,citecolor=BrickRed,linkcolor=NavyBlue,urlcolor=RoyalBlue]{hyperref}
\usepackage{orcidlink}

\numberwithin{equation}{section}

\newcommand{\keywords}[1]{\par\medskip\noindent\textbf{Keywords.} #1}
\newcommand{\subjclass}[1]{\par\smallskip\noindent\textbf{Mathematics Subject Classification (2020).} #1}

\usepackage{algorithm}

\usepackage[export]{adjustbox} 

\makeatletter
\newenvironment{breakablealgorithm}
  {
   \begin{center}
     \refstepcounter{algorithm}
     \hrule height.8pt depth0pt \kern2pt
     \renewcommand{\caption}[2][\relax]{
       {\raggedright\textbf{\fname@algorithm~\thealgorithm} ##2\par}%
       \ifx\relax##1\relax 
         \addcontentsline{loa}{algorithm}{\protect\numberline{\thealgorithm}##2}%
       \else 
         \addcontentsline{loa}{algorithm}{\protect\numberline{\thealgorithm}##1}%
       \fi
       \kern2pt\hrule\kern2pt
     }
  }{
     \kern2pt\hrule\relax
   \end{center}
  }
\makeatother

\usepackage{algorithmic}

\title{Derivative-Free Structured Updates for Muon}
\newcommand{\name}[1]{#1}
\author{
\name{Pengcheng Xie\textsuperscript{a}\thanks{CONTACT P. Xie. Email: \href{mailto:pxie@lbl.gov}{pxie@lbl.gov}}\,\orcidlink{0000-0001-5973-1535}}\\[0.5ex]
\small\textsuperscript{a}Lawrence Berkeley National Laboratory \& University of California\\
\small Berkeley, CA 94720, USA
}
\date{}
\hypersetup{pdftitle={Derivative-Free Structured Updates for Muon},
  pdfauthor={Pengcheng Xie},
  pdfsubject={Derivative-free optimization of matrix-valued parameters}}

\begin{document}

\maketitle

\begin{abstract}
Muon updates matrix-valued neural-network parameters by orthogonalizing a
gradient-based momentum matrix. Its reliance on derivatives limits its use when
gradients are unavailable or unreliable. We develop a derivative-free framework
that constructs Muon-style updates from structured finite differences. Four
variants are considered: full entrywise recovery, random low-rank surrogates,
basis-aligned rank-one probing, and direct structured search. Exhaustive
basis-aligned probing is equivalent, up to positive scaling before ideal polar
orthogonalization, to coordinate finite differences. Matrix-regression
experiments show that random rank-one probing can reduce the number of function
evaluations substantially, at the cost of less accurate updates. Controlled
noisy-gradient experiments on regression and a neural network illustrate when
accurate function values can compensate for an unreliable gradient oracle.
A small CartPole study further examines orthogonal rank-one probes under a
fixed episode budget. These results support structured probing as a practical
option for selected black-box problems; they do not establish a general
convergence guarantee or an advantage over accurate, inexpensive gradients.
\end{abstract}

\section{Introduction}

Muon is an optimizer for matrix-valued neural-network parameters
\cite{jordan2024muon}. It forms momentum from the gradient and orthogonalizes
the resulting matrix before updating the parameters. We study how this update
can be constructed using function evaluations when derivatives are unavailable
or unreliable.

The orthogonalized update has been interpreted as a layerwise duality map
\cite{bernstein2024modular} and as a steepest-descent direction in the spectral
norm \cite{li2025note}. Convergence analyses identify favorable low-rank and
block-structured curvature regimes \cite{shen2025convergence}, while decoupled
weight decay has been related to a spectral-norm constraint
\cite{chen2025spectral}. Muon has also been used in multi-billion-parameter
language-model training \cite{liu2025scalable}.

Su's isotropic curvature model provides a complementary analysis of the
singular values of an update matrix \cite{su2025isotropic}. Under its curvature
growth assumptions, the optimal update has a more homogeneous spectrum than
the gradient; full orthogonalization is optimal under a particular
phase-transition condition. Du and Su derive Newton-Muon from a quadratic
surrogate, with right preconditioning by the inverse input second-moment
matrix before orthogonalization \cite{du2026newtonmuon}. These analyses show
how curvature assumptions affect the choice of matrix update. Our focus is on
constructing Muon-style directions from function values when gradient
information is unavailable or unreliable.

At each iteration, Muon forms a stochastic gradient matrix $\boldsymbol{G}_t$, incorporates
momentum to obtain $\boldsymbol{M}_t$, and applies a short Newton--Schulz iteration
\cite{higham2008functions} to produce an approximately semi-orthogonal matrix
$\boldsymbol{O}_t$. The parameter update has the form
\[
\boldsymbol{W}_{t+1} = \boldsymbol{W}_t - \eta \boldsymbol{O}_t.
\]
For a compact singular value decomposition
$\boldsymbol{M}_t=\boldsymbol{U}_t\boldsymbol{\Sigma}_t\boldsymbol{V}_t^\top$,
the idealized update is $\boldsymbol{U}_t\boldsymbol{V}_t^\top$ on the nonzero
singular subspaces. This transformation retains paired singular directions
while removing their singular-value magnitudes. It motivates seeking useful
update directions from partial matrix information. However, a small error in a
low-rank gradient approximation need not imply a small error in the polar
factor: even small nonzero singular values contribute unit-weight directions.
The usefulness of low-rank probing must therefore be assessed through the
resulting optimization behavior.

When an objective is evaluated by a simulator or an external solver, its
derivatives may be inaccessible. Derivative-free methods use function values
in their place \cite{conn2009introduction,spall1992multivariate}.

Related model-based methods reduce the information needed to construct local
approximations. Xie and Yuan study underdetermined quadratic interpolation
models that account for trust-region iteration properties
\cite{xie2025underdetermined}, and propose least $H^2$ norm updates that control
the change between successive quadratic models \cite{xie2026h2}. Their
2D-MoSub method uses two-dimensional model-based subspaces for large-scale
derivative-free optimization \cite{xie2026mosub}. Here we use function values to approximate matrix update directions.

Finite differences give a direct way to approximate the gradient matrix.
Accurate full reconstruction generally requires a measurement budget that grows
with the number of parameters. Randomized estimators use fewer queries but
introduce sampling error \cite{nesterov2017random}. We use rank-one
perturbations to construct low-rank gradient surrogates and compare the induced
Muon updates with exhaustive coordinate probing. We also consider direct
search over rank-one directions. These methods trade query cost against
update quality; the numerical experiments measure both.

Section~2 describes the algorithms. Section~3 reports results for matrix
regression, neural-network training, and CartPole policy optimization.
Section~4 discusses the conclusions and remaining questions.

\begin{breakablealgorithm}
\caption{Muon}
\begin{algorithmic}\label{alg:01}
\REQUIRE{$\eta > 0, \lambda>0, \beta \in [0,1), \boldsymbol{M}_{-1} := \boldsymbol{0}, \boldsymbol{W}_0 \in \mathbb{R}^{m\times n}$}
\FOR{$t=0$ to $T-1$}
\STATE{$\boldsymbol{M}_t := \beta \boldsymbol{M}_{t-1} + (1-\beta) \nabla f_{\mathcal{S}_t}(\boldsymbol{W}_t)$}
\IF{(Nesterov = True)}
\STATE{$\boldsymbol{C}_t := \beta \boldsymbol{M}_t + (1-\beta) \nabla f_{\mathcal{S}_t}(\boldsymbol{W}_t)$}
\ELSE
\STATE{$\boldsymbol{C}_t := \boldsymbol{M}_t$}
\ENDIF
\STATE{$\boldsymbol{O}_t := \text{NewtonSchulz}(\boldsymbol{C}_t)$}
\IF{(weight decay = True)}
\STATE{$\boldsymbol{W}_{t+1} := (1-\eta \lambda) \boldsymbol{W}_t - \eta \boldsymbol{O}_t$}
\ELSE
\STATE{$\boldsymbol{W}_{t+1} := \boldsymbol{W}_t - \eta \boldsymbol{O}_t$}
\ENDIF
\ENDFOR \\
\RETURN{$\boldsymbol{W}_T$}
\end{algorithmic}
\end{breakablealgorithm}

\subsection{Rank-One Measurements}

For $\boldsymbol{W}\in\mathbb{R}^{m\times n}$, coordinate finite differences
require $mn$ perturbed function values. If the gradient is well approximated
by a rank-$r$ matrix, its factors contain only $r(m+n)$ entries. This suggests
using measurements adapted to the matrix structure, although it does not by
itself guarantee accurate recovery of the polar factor.

For a rank-one perturbation $\boldsymbol{H}=\boldsymbol{u}\boldsymbol{v}^\top$,
\[
\frac{f(\boldsymbol{W}+h\boldsymbol{u}\boldsymbol{v}^\top)-f(\boldsymbol{W})}{h}
\approx
\langle\nabla f(\boldsymbol{W}),\boldsymbol{u}\boldsymbol{v}^\top\rangle_F
=
\boldsymbol{u}^\top\nabla f(\boldsymbol{W})\boldsymbol{v}.
\]
Once the base value is cached, each measurement costs one additional function
evaluation. The factors can be drawn randomly, chosen from coordinate bases,
or adapted to previous measurements. The methods below use these measurements
to assemble a momentum surrogate or select a trial direction.

The query cost is separate from the cost of orthogonalization. Replacing a
gradient by a small set of probes reduces the number of measurements, but the
subsequent Newton--Schulz iteration still requires matrix multiplications.
Low-rank storage could reduce that cost if the momentum is kept in factored
form; the methods tested here store the momentum as a full matrix.

\subsection{A Concrete Sparse Matrix Example}

A masked quadratic objective gives a simple example of the possible query savings.
Consider the optimization variable \(\boldsymbol{W}\in\mathbb{R}^{4\times 4}\) and the
black-box objective
\begin{equation}
\label{eq:sparse_example_objective}
 f(\boldsymbol{W})
 =\frac{1}{2}(\boldsymbol{W}_{12}-3)^2+\frac{1}{2}(\boldsymbol{W}_{34}+1)^2.
\end{equation}
This is a masked matrix-fitting problem: although \(\boldsymbol{W}\) has \(16\) entries, the
loss observes only the sparse set

\[
 \Omega=\{(1,2),(3,4)\}.
\]

Such a mask can represent, for example, a matrix completion problem with only
a few observed entries, a sparse interaction model, or a simulator in which
only selected input--output channels are coupled.  Starting from \(\boldsymbol{W}_0=0\), the
objective and gradient are
\[
 f(\boldsymbol{W}_0)=\frac{1}{2}(3^2+1^2)=5,
 \qquad
 \boldsymbol{G}_0=\nabla f(\boldsymbol{W}_0)=-3\boldsymbol{E}_{12}+\boldsymbol{E}_{34}
 =
 \begin{bmatrix}
 0&-3&0&0\\
 0& 0&0&0\\
 0& 0&0&1\\
 0& 0&0&0
 \end{bmatrix},
\]
where \(\boldsymbol{E}_{ij}=\boldsymbol{e}_i \boldsymbol{e}_j^\top\).  Thus only two of the \(16\) gradient entries
contain any first-order information.

For clarity, take no momentum at the first iteration and replace the finite
Newton--Schulz approximation by the exact polar factor.  A compact singular
value decomposition of \(\boldsymbol{G}_0\) gives singular values \(3\) and \(1\), and hence
the ideal Muon direction is
\begin{equation}
\label{eq:sparse_example_muon_direction}
 \boldsymbol{O}_0=\operatorname{polar}(\boldsymbol{G}_0)=-\boldsymbol{E}_{12}+\boldsymbol{E}_{34}
 =
 \begin{bmatrix}
 0&-1&0&0\\
 0& 0&0&0\\
 0& 0&0&1\\
 0& 0&0&0
 \end{bmatrix}.
\end{equation}
Muon discards the unequal magnitudes \(3\) and \(1\), but retains their signs
and the two active row--column pairs.  With step size \(\eta=1\),
\[
 \boldsymbol{W}_1=\boldsymbol{W}_0-\boldsymbol{O}_0=\boldsymbol{E}_{12}-\boldsymbol{E}_{34},
 \qquad
 f(\boldsymbol{W}_1)=\frac{1}{2}(1-3)^2+\frac{1}{2}(-1+1)^2=2.
\]
The objective therefore decreases from \(5\) to \(2\), and all \(14\) inactive
coordinates play no role in this Muon step.

The derivative-free implication can also be quantified exactly.  A naive
central-difference reconstruction uses
\[
 \widehat{\boldsymbol{G}}_{ij}
 =\frac{f(\boldsymbol{W}_0+h\boldsymbol{E}_{ij})-f(\boldsymbol{W}_0-h\boldsymbol{E}_{ij})}{2h}
\]
for every \(i,j\), requiring \(2mn=32\) function evaluations.  In contrast, if
the mask \(\Omega\) is known from the problem structure, only the four values
associated with \(\boldsymbol{E}_{12}\) and \(\boldsymbol{E}_{34}\) are needed.  Because
\eqref{eq:sparse_example_objective} is quadratic, central differences are exact in exact arithmetic
for every \(h>0\):
\[
 \frac{f(h\boldsymbol{E}_{12})-f(-h\boldsymbol{E}_{12})}{2h}=-3,
 \qquad
 \frac{f(h\boldsymbol{E}_{34})-f(-h\boldsymbol{E}_{34})}{2h}=1.
\]
These two measurements already determine \eqref{eq:sparse_example_muon_direction};
recovering the other \(14\) zero entries provides no additional information to
the update.  More generally, if only \(s\ll mn\) entries can be active and their
support is known, the query count falls from \(2mn\) to \(2s\).  If the support
is unknown, random signed matrix probes
\[
 y_k=\frac{f(\boldsymbol{W}+h\boldsymbol{H}_k)-f(\boldsymbol{W}-h\boldsymbol{H}_k)}{2h}
 \approx \langle \nabla f(\boldsymbol{W}),\boldsymbol{H}_k\rangle_F
\]
turn the task into sparse matrix sensing: the informative support and signs can
be estimated from a number of aggregate measurements governed by the sparsity
level rather than by blindly probing every coordinate.

Here the active coordinates and their signs determine the Muon direction.
The reduction from 32 to four queries uses prior knowledge of the mask.
With unknown support, additional measurements and a recovery procedure are
needed; sparsity alone does not provide the support.

\section{Algorithms for Derivative-Free Structured Muon}

This section presents four algorithmic realizations of the proposed
derivative-free structured Muon framework. The methods differ in how much
gradient information they attempt to recover. The first reconstructs an
entrywise approximation of the full gradient matrix and serves as a baseline.
The second estimates a low-rank structured approximation using randomly
or otherwise globally chosen rank-one probes. The third, denoted V2, retains
the rank-one factorization but aligns the left and right probing vectors with
the coordinate bases. The fourth bypasses gradient estimation and searches
directly for structured update directions.

\subsection{Entrywise-Recovery Derivative-Free Muon}

The most direct derivative-free extension of Muon approximates the full
gradient matrix by finite differences. Let \(\boldsymbol{W}_t\in\mathbb{R}^{m\times n}\) be
the current matrix parameter, and let \(\boldsymbol{H}_i\in\mathbb{R}^{m\times n}\) be
probing matrices. For each \(\boldsymbol{H}_i\), compute the directional finite-difference
measurement
\[
d_{t,i}
=
\frac{f(\boldsymbol{W}_t+h\boldsymbol{H}_i)-f(\boldsymbol{W}_t)}{h}
\approx
\langle \nabla f(\boldsymbol{W}_t),\boldsymbol{H}_i\rangle_F .
\]
The full gradient matrix is then reconstructed through a least-squares problem,
and the resulting approximation is passed to the standard Muon
orthogonalization step.

\begin{breakablealgorithm}
\caption{Entrywise-Recovery Derivative-Free Muon}
\label{alg:df_muon_entrywise}
\begin{algorithmic}
\REQUIRE{Step size \(\eta>0\), weight decay parameter \(\lambda\ge 0\), momentum parameter \(\beta\in[0,1)\), finite-difference radius \(h>0\), probing matrices \(\{\boldsymbol{H}_i\}_{i=1}^p\subset\mathbb{R}^{m\times n}\), initial matrix \(\boldsymbol{W}_0\in\mathbb{R}^{m\times n}\), \(\boldsymbol{M}_{-1}=0\)}
\FOR{\(t=0,1,\ldots,T-1\)}
    \STATE Evaluate \(f(\boldsymbol{W}_t)\).
    \FOR{\(i=1,\ldots,p\)}
        \STATE Evaluate \(f(\boldsymbol{W}_t+h\boldsymbol{H}_i)\).
        \STATE Compute
        \[
        d_{t,i}:=\frac{f(\boldsymbol{W}_t+h\boldsymbol{H}_i)-f(\boldsymbol{W}_t)}{h}.
        \]
    \ENDFOR
    \STATE Compute the full matrix estimate
    \[
    \widehat{\boldsymbol{G}}_t
    \in
    \arg\min_{\boldsymbol{X}\in\mathbb{R}^{m\times n}}
    \sum_{i=1}^p
    \bigl(\langle \boldsymbol{X},\boldsymbol{H}_i\rangle_F-d_{t,i}\bigr)^2 .
    \]
    \STATE Update the momentum matrix
    \[
    \boldsymbol{M}_t:=\beta \boldsymbol{M}_{t-1}+(1-\beta)\widehat{\boldsymbol{G}}_t .
    \]
    \STATE Compute the structured Muon direction
    \[
    \boldsymbol{O}_t:=\operatorname{NewtonSchulz}(\boldsymbol{M}_t).
    \]
    \IF{weight decay is used}
        \STATE
        \[
        \boldsymbol{W}_{t+1}:=(1-\eta\lambda)\boldsymbol{W}_t-\eta \boldsymbol{O}_t .
        \]
    \ELSE
        \STATE
        \[
        \boldsymbol{W}_{t+1}:=\boldsymbol{W}_t-\eta \boldsymbol{O}_t .
        \]
    \ENDIF
\ENDFOR
\RETURN{\(\boldsymbol{W}_T\)}
\end{algorithmic}
\end{breakablealgorithm}

This method is conceptually simple but expensive. In general, accurate
entrywise reconstruction of a matrix in \(\mathbb{R}^{m\times n}\) requires a
number of probing directions comparable to \(mn\). The structured methods below
are designed to avoid this cost.

\subsection{Low-Rank Structured Derivative-Free Muon}

The entrywise-recovery method approximates the entire gradient matrix, whereas
Muon ultimately uses the orthogonalized structure of the momentum matrix rather
than the exact entrywise values of the gradient. This motivates a low-rank
derivative-free approximation, in the spirit of randomized low-rank matrix
approximation methods \cite{halko2011finding}.

We consider rank-one probing directions
\[
\boldsymbol{H}_i=\boldsymbol{u}_i \boldsymbol{v}_i^\top,
\qquad
\boldsymbol{u}_i\in\mathbb{R}^m,\quad \boldsymbol{v}_i\in\mathbb{R}^n.
\]
The corresponding finite-difference measurement satisfies
\[
\frac{f(\boldsymbol{W}_t+h \boldsymbol{u}_i \boldsymbol{v}_i^\top)-f(\boldsymbol{W}_t)}{h}
\approx
\langle \nabla f(\boldsymbol{W}_t),\boldsymbol{u}_i \boldsymbol{v}_i^\top\rangle_F
=
\boldsymbol{u}_i^\top \nabla f(\boldsymbol{W}_t)\boldsymbol{v}_i .
\]
With a cached base value, each additional function evaluation provides an
approximate bilinear measurement of the gradient matrix.

If the vectors \(\boldsymbol{u}_i\) and \(\boldsymbol{v}_i\) are normalized, then a simple low-rank
surrogate can be formed as
\[
\widehat{\boldsymbol{G}}_t
=
\sum_{i=1}^r
\alpha_{t,i} \boldsymbol{u}_i \boldsymbol{v}_i^\top,
\qquad
\alpha_{t,i}
=
\frac{f(\boldsymbol{W}_t+h \boldsymbol{u}_i \boldsymbol{v}_i^\top)-f(\boldsymbol{W}_t)}{h}.
\]
If \(\boldsymbol{u}_i\) and \(\boldsymbol{v}_i\) are not normalized, one may instead use
\[
\alpha_{t,i}
=
\frac{f(\boldsymbol{W}_t+h \boldsymbol{u}_i \boldsymbol{v}_i^\top)-f(\boldsymbol{W}_t)}
{h\|\boldsymbol{u}_i\|_2^2\|\boldsymbol{v}_i\|_2^2}.
\]

This sum is a rank-at-most-$r$ surrogate, not generally the orthogonal
projection onto the span of the probes unless the atoms are mutually
orthogonal. Nor do $r$ arbitrary scalar measurements identify a general
rank-$r$ matrix. The construction is used as a heuristic update surrogate.
Although each surrogate has low rank, accumulating momentum over different
probes can increase the rank of the momentum matrix.

\begin{breakablealgorithm}
\caption{Low-Rank Structured Derivative-Free Muon}
\label{alg:df_muon_lowrank}
\begin{algorithmic}
\REQUIRE{Step size \(\eta>0\), weight decay parameter \(\lambda\ge 0\), momentum parameter \(\beta\in[0,1)\), finite-difference radius \(h>0\), target rank \(r\), probing pairs \(\{(\boldsymbol{u}_i,\boldsymbol{v}_i)\}_{i=1}^r\), initial matrix \(\boldsymbol{W}_0\in\mathbb{R}^{m\times n}\), \(\boldsymbol{M}_{-1}=0\)}
\FOR{\(t=0,1,\ldots,T-1\)}
    \STATE Evaluate \(f(\boldsymbol{W}_t)\).
    \FOR{\(i=1,\ldots,r\)}
        \STATE Form the rank-one probing matrix
        \[
        \boldsymbol{H}_i:=\boldsymbol{u}_i \boldsymbol{v}_i^\top .
        \]
        \STATE Evaluate \(f(\boldsymbol{W}_t+h\boldsymbol{H}_i)\).
        \STATE Compute
        \[
        \alpha_{t,i}
        :=
        \frac{f(\boldsymbol{W}_t+h\boldsymbol{H}_i)-f(\boldsymbol{W}_t)}
        {h\|\boldsymbol{H}_i\|_F^2}.
        \]
    \ENDFOR
    \STATE Construct the low-rank gradient surrogate
    \[
    \widehat{\boldsymbol{G}}_t
    :=
    \sum_{i=1}^r \alpha_{t,i} \boldsymbol{H}_i .
    \]
    \STATE Update the momentum matrix
    \[
    \boldsymbol{M}_t:=\beta \boldsymbol{M}_{t-1}+(1-\beta)\widehat{\boldsymbol{G}}_t .
    \]
    \STATE Compute the Muon-style structured direction
    \[
    \boldsymbol{O}_t:=\operatorname{NewtonSchulz}(\boldsymbol{M}_t).
    \]
    \IF{weight decay is used}
        \STATE
        \[
        \boldsymbol{W}_{t+1}:=(1-\eta\lambda)\boldsymbol{W}_t-\eta \boldsymbol{O}_t .
        \]
    \ELSE
        \STATE
        \[
        \boldsymbol{W}_{t+1}:=\boldsymbol{W}_t-\eta \boldsymbol{O}_t .
        \]
    \ENDIF
\ENDFOR
\RETURN{\(\boldsymbol{W}_T\)}
\end{algorithmic}
\end{breakablealgorithm}

\subsection{Basis-Aligned Rank-One DFO-Muon (V2)}
\label{subsec:df_muon_v2}

The basis-aligned variant, denoted V2 in the experiments, chooses both probing
factors from the canonical bases. Let

\[
\{\boldsymbol{e}^{(m)}_i\}_{i=1}^m
\quad\text{and}\quad
\{\boldsymbol{e}^{(n)}_j\}_{j=1}^n
\]
denote the standard bases of \(\mathbb{R}^m\) and \(\mathbb{R}^n\), respectively,
and set
\[
\boldsymbol{u}_{ij}=\boldsymbol{e}^{(m)}_i,
\qquad
\boldsymbol{v}_{ij}=\boldsymbol{e}^{(n)}_j,
\qquad
\boldsymbol{H}_{ij}
=
\boldsymbol{u}_{ij}\boldsymbol{v}_{ij}^{\top}
=
\boldsymbol{e}^{(m)}_i(\boldsymbol{e}^{(n)}_j)^{\top}.
\]
Each atom has rank one and probes a single matrix entry.

For a coordinate pair \((i,j)\), V2 uses the rank-one perturbation
\[
\boldsymbol{D}_{ij}=h\boldsymbol{H}_{ij}
=h\boldsymbol{e}^{(m)}_i(\boldsymbol{e}^{(n)}_j)^{\top}
\]
and the one-sided directional measurement
\begin{equation}
\alpha_{t,ij}
=
\frac{f(\boldsymbol{W}_t+\boldsymbol{D}_{ij})
      -f(\boldsymbol{W}_t)}{h}
=
\frac{f(\boldsymbol{W}_t+h\boldsymbol{H}_{ij})
      -f(\boldsymbol{W}_t)}{h}.
\label{eq:v2_alpha}
\end{equation}
Since \(\|\boldsymbol{H}_{ij}\|_F=1\), Taylor expansion gives
\[
\alpha_{t,ij}
=
\left\langle \nabla f(\boldsymbol{W}_t),
\boldsymbol{H}_{ij}\right\rangle_F+O(h)
=
\bigl[\nabla f(\boldsymbol{W}_t)\bigr]_{ij}+O(h).
\]
For the regularized matrix-regression objective used in Section~3, this
relation can be written exactly.  Define
\[
\boldsymbol{A}:=\frac{\boldsymbol{X}\boldsymbol{X}^{\top}}{N}.
\]
Because the objective is quadratic and
\(\nabla f(\boldsymbol{W})=(\boldsymbol{W}\boldsymbol{X}-\boldsymbol{Y})
\boldsymbol{X}^{\top}/N+\lambda\boldsymbol{W}\), one obtains
\begin{equation}
\alpha_{t,ij}
=
\bigl[\nabla f(\boldsymbol{W}_t)\bigr]_{ij}
+\frac{h}{2}\bigl(A_{jj}+\lambda\bigr).
\label{eq:v2_exact_quadratic_fd}
\end{equation}
Thus, the one-sided V2 measurement has an explicit \(O(h)\) bias; for a fixed
column \(j\), the bias is the same for every row \(i\).  This bias comes from
the finite difference itself and is separate from the rank-one assembly rule
described next.

V2 accumulates the scaled perturbations to form
\begin{equation}
\widetilde{\boldsymbol{G}}^{\mathrm{V2}}_t
=
c_{m,n}
\sum_{i=1}^{m}\sum_{j=1}^{n}
\alpha_{t,ij}\boldsymbol{D}_{ij},
\qquad
c_{m,n}=\frac{mn}{\operatorname{size}(\boldsymbol{W}_t)}.
\label{eq:v2_surrogate}
\end{equation}
For an \(m\)-by-\(n\) matrix,
\(\operatorname{size}(\boldsymbol{W}_t)=mn\), so \(c_{m,n}=1\).  Consequently,
if
\[
\widehat{\boldsymbol{G}}^{\mathrm{coord}}_t
=
\sum_{i=1}^{m}\sum_{j=1}^{n}
\alpha_{t,ij}\boldsymbol{H}_{ij}
\]
denotes the usual coordinate finite-difference matrix, then the surrogate used
by V2 satisfies the exact identity
\begin{equation}
\widetilde{\boldsymbol{G}}^{\mathrm{V2}}_t
=h\widehat{\boldsymbol{G}}^{\mathrm{coord}}_t.
\label{eq:v2_scaling_identity}
\end{equation}
Thus, exhaustive V2 uses the same measurements as coordinate finite
differences and rescales the reconstructed matrix by $h$. Both methods require
$mn$ perturbed evaluations per iteration.

The V2 momentum and parameter updates are
\begin{align}
\boldsymbol{M}^{\mathrm{V2}}_t
&=
\beta\boldsymbol{M}^{\mathrm{V2}}_{t-1}
+(1-\beta)\widetilde{\boldsymbol{G}}^{\mathrm{V2}}_t,
\label{eq:v2_momentum}\\
\boldsymbol{O}^{\mathrm{V2}}_t
&=
\operatorname{Polar}(\boldsymbol{M}^{\mathrm{V2}}_t),
\qquad
\boldsymbol{W}_{t+1}
=
\boldsymbol{W}_t-\eta\boldsymbol{O}^{\mathrm{V2}}_t.
\label{eq:v2_update}
\end{align}
Here \(\operatorname{Polar}\) may be evaluated by an SVD or approximated by a
Newton--Schulz iteration.  If V2 and coordinate recovery start from the same
zero momentum and follow the same parameter trajectory, then
\[
\boldsymbol{M}^{\mathrm{V2}}_t
=h\boldsymbol{M}^{\mathrm{coord}}_t.
\]
The ideal polar map is homogeneous of degree zero for positive scalars:
\[
\operatorname{Polar}(h\boldsymbol{M})
=\operatorname{Polar}(\boldsymbol{M}),
\qquad h>0.
\]
Under exhaustive probing and consistent polar-factor conventions, the two
methods therefore give the same direction in exact arithmetic. Finite-precision
trajectories can differ, particularly near rank deficiency, where the polar
factor is sensitive to perturbations. The absolute zero threshold used by the
implementation can also respond differently to the scaled momentum.

\begin{breakablealgorithm}
\caption{Basis-Aligned Rank-One DFO-Muon (V2)}
\label{alg:df_muon_v2}
\begin{algorithmic}
\REQUIRE{Step size \(\eta>0\), momentum parameter \(\beta\in[0,1)\), finite-difference radius \(h>0\), initial matrix \(\boldsymbol{W}_0\in\mathbb{R}^{m\times n}\), \(\boldsymbol{M}_{-1}=0\)}
\FOR{\(t=0,1,\ldots,T-1\)}
    \STATE Evaluate \(b_t:=f(\boldsymbol{W}_t)\).
    \STATE Set \(\widetilde{\boldsymbol{G}}^{\mathrm{V2}}_t:=\boldsymbol{0}\).
    \FOR{\(i=1,\ldots,m\)}
        \FOR{\(j=1,\ldots,n\)}
            \STATE Choose
            \(
            \boldsymbol{u}_{ij}:=\boldsymbol{e}^{(m)}_i
            \) and
            \(
            \boldsymbol{v}_{ij}:=\boldsymbol{e}^{(n)}_j
            \).
            \STATE Form the scaled rank-one perturbation
            \[
            \boldsymbol{D}_{ij}
            :=h\boldsymbol{u}_{ij}\boldsymbol{v}_{ij}^{\top}.
            \]
            \STATE Evaluate \(f(\boldsymbol{W}_t+\boldsymbol{D}_{ij})\) and compute
            \[
            \alpha_{t,ij}
            :=
            \frac{f(\boldsymbol{W}_t+\boldsymbol{D}_{ij})-b_t}{h}.
            \]
            \STATE Accumulate the rank-one atom
            \[
            \widetilde{\boldsymbol{G}}^{\mathrm{V2}}_t
            :=
            \widetilde{\boldsymbol{G}}^{\mathrm{V2}}_t
            +\alpha_{t,ij}\boldsymbol{D}_{ij}.
            \]
        \ENDFOR
    \ENDFOR
    \STATE Apply the scaling
    \[
    \widetilde{\boldsymbol{G}}^{\mathrm{V2}}_t
    :=
    \frac{mn}{\operatorname{size}(\boldsymbol{W}_t)}
    \widetilde{\boldsymbol{G}}^{\mathrm{V2}}_t.
    \]
    \STATE Update the momentum
    \[
    \boldsymbol{M}_t
    :=
    \beta\boldsymbol{M}_{t-1}
    +(1-\beta)\widetilde{\boldsymbol{G}}^{\mathrm{V2}}_t.
    \]
    \STATE Compute
    \(
    \boldsymbol{O}_t:=\operatorname{Polar}(\boldsymbol{M}_t)
    \)
    and update
    \[
    \boldsymbol{W}_{t+1}:=\boldsymbol{W}_t-\eta\boldsymbol{O}_t.
    \]
\ENDFOR
\RETURN{\(\boldsymbol{W}_T\)}
\end{algorithmic}
\end{breakablealgorithm}

Each V2 probe is rank one, but the sum in
\eqref{eq:v2_surrogate} can have rank as large as \(\min\{m,n\}\).  Therefore,
``low rank'' in the name refers to the factorized rank-one measurement atoms,
not to a guarantee that the final surrogate has rank one or a prescribed small
rank.  The exhaustive implementation uses \(mn\) perturbed evaluations plus
one base evaluation per iteration, giving
\[
1+T(mn+1)
\]
recorded objective evaluations over \(T\) iterations when the initial
evaluation is included.  Hence its present query complexity matches exhaustive
coordinate probing rather than the \(r\)-probe randomized low-rank method.

A reduced-query variant would replace the full Cartesian set
\(\{1,\ldots,m\}\times\{1,\ldots,n\}\) by a subset
\(\mathcal{I}_t\) of \(q\ll mn\) coordinate pairs, selected randomly,
cyclically, or adaptively.  For uniform random coordinate sampling, a surrogate
of the form
\[
\frac{mn}{q}
\sum_{(i,j)\in\mathcal{I}_t}
\alpha_{t,ij}\boldsymbol{H}_{ij}
\]
gives an estimator of the coordinate finite-difference matrix. This sampled
variant is not included in the experiments below.

\subsection{Direct Structured Update Search}

We next consider a derivative-free method that does not construct a gradient
approximation. Instead, it searches directly over structured update directions.
At iteration \(t\), generate candidate rank-one directions
\[
\boldsymbol{H}_i=\boldsymbol{u}_i \boldsymbol{v}_i^\top,
\qquad i=1,\ldots,q,
\]
and evaluate their objective decrease through
\[
s_{t,i}
=
\frac{f(\boldsymbol{W}_t+h\boldsymbol{H}_i)-f(\boldsymbol{W}_t)}{h}.
\]
A smaller value of \(s_{t,i}\) indicates a better descent direction, so the
selected structured direction is
\[
\boldsymbol{H}_t
=
\boldsymbol{H}_{i_t},
\qquad
i_t\in\arg\min_{1\le i\le q}s_{t,i}.
\]
Using the candidate with the smallest score, the
parameter update is
\[
\boldsymbol{W}_{t+1}=\boldsymbol{W}_t+\eta \boldsymbol{H}_t.
\]

\begin{breakablealgorithm}
\caption{Direct Structured Derivative-Free Update Search}
\label{alg:df_direct_structured}
\begin{algorithmic}
\REQUIRE{Step size \(\eta>0\), finite-difference radius \(h>0\), number of candidates \(q\), initial matrix \(\boldsymbol{W}_0\in\mathbb{R}^{m\times n}\)}
\FOR{\(t=0,1,\ldots,T-1\)}
    \STATE Evaluate \(f(\boldsymbol{W}_t)\).
    \FOR{\(i=1,\ldots,q\)}
        \STATE Sample candidate vectors \(\boldsymbol{u}_i\in\mathbb{R}^m\) and \(\boldsymbol{v}_i\in\mathbb{R}^n\).
        \STATE Form the rank-one direction
        \[
        \boldsymbol{H}_i:=\boldsymbol{u}_i \boldsymbol{v}_i^\top .
        \]
        \STATE Evaluate \(f(\boldsymbol{W}_t+h\boldsymbol{H}_i)\).
        \STATE Compute the directional score
        \[
        s_{t,i}
        :=
        \frac{f(\boldsymbol{W}_t+h\boldsymbol{H}_i)-f(\boldsymbol{W}_t)}{h}.
        \]
    \ENDFOR
    \STATE Select the best structured direction
    \[
    i_t\in\arg\min_{1\le i\le q}s_{t,i},
    \qquad
    \boldsymbol{H}_t:=\boldsymbol{H}_{i_t}.
    \]
    \STATE Update
    \[
    \boldsymbol{W}_{t+1}:=\boldsymbol{W}_t+\eta \boldsymbol{H}_t .
    \]
\ENDFOR
\RETURN{\(\boldsymbol{W}_T\)}
\end{algorithmic}
\end{breakablealgorithm}

The always-move rule does not guarantee descent: all candidate scores may be
positive, and the update radius $\eta$ may differ from the probing radius $h$.
An acceptance test or a line search would be needed to enforce actual decrease.

The method can also be modified to use a symmetric two-point search. In that
case, for each candidate \(\boldsymbol{H}_i\), evaluate both \(f(\boldsymbol{W}_t+h\boldsymbol{H}_i)\) and
\(f(\boldsymbol{W}_t-h\boldsymbol{H}_i)\), and then choose between \(\boldsymbol{H}_i\) and \(-\boldsymbol{H}_i\). This yields the
score
\[
s_{t,i}^{\pm}
=
\min\left\{
f(\boldsymbol{W}_t+h\boldsymbol{H}_i),
f(\boldsymbol{W}_t-h\boldsymbol{H}_i)
\right\}.
\]
The symmetric variant is often more robust because it does not fix the sign of
the sampled rank-one perturbation in advance.

\section{Numerical Experiments}

We now present a concrete matrix optimization problem arising from the training
of a single-layer neural network. The example illustrates how derivative-free
structured Muon methods can be applied to matrix-valued parameters.

Let
\[
\boldsymbol{x}_j\in\mathbb{R}^n,\qquad \boldsymbol{y}_j\in\mathbb{R}^m,
\qquad j=1,\ldots,N,
\]
be the training data, where \(\boldsymbol{x}_j\) is an input vector and \(\boldsymbol{y}_j\) is the
corresponding target output. Consider a single-layer neural network of the
form
\[
\Phi_{\boldsymbol{W}}(\boldsymbol{x})
=\sigma(\boldsymbol{W}\boldsymbol{x}),
\]
where
\[
\boldsymbol{W}\in\mathbb{R}^{m\times n}
\]
is the trainable weight matrix and \(\sigma:\mathbb{R}^m\to\mathbb{R}^m\) is an
activation function applied componentwise. The empirical risk minimization
problem is
\[
\min_{\boldsymbol{W}\in\mathbb{R}^{m\times n}}
f(\boldsymbol{W})
:=
\frac{1}{2N}
\sum_{j=1}^N
\|\sigma(\boldsymbol{W}\boldsymbol{x}_j)-\boldsymbol{y}_j\|_2^2
+
\frac{\lambda}{2}\|\boldsymbol{W}\|_F^2 .
\]
The optimization variable is therefore the weight matrix \(\boldsymbol{W}\).

For example, if \(\sigma\) is the identity map, then the problem reduces to
regularized linear regression:
\[
f(\boldsymbol{W})
=
\frac{1}{2N}
\sum_{j=1}^N
\|\boldsymbol{W}\boldsymbol{x}_j-\boldsymbol{y}_j\|_2^2
+
\frac{\lambda}{2}\|\boldsymbol{W}\|_F^2 .
\]
Let
\[
\boldsymbol{X}=
\begin{bmatrix}
\boldsymbol{x}_1 & \boldsymbol{x}_2 & \cdots & \boldsymbol{x}_N
\end{bmatrix}
\in\mathbb{R}^{n\times N},
\qquad
\boldsymbol{Y}=
\begin{bmatrix}
\boldsymbol{y}_1 & \boldsymbol{y}_2 & \cdots & \boldsymbol{y}_N
\end{bmatrix}
\in\mathbb{R}^{m\times N}.
\]
Then the objective can be written compactly as
\[
f(\boldsymbol{W})
=
\frac{1}{2N}\|\boldsymbol{W}\boldsymbol{X}-\boldsymbol{Y}\|_F^2
+
\frac{\lambda}{2}\|\boldsymbol{W}\|_F^2 .
\]
In this special case, the exact gradient is
\[
\nabla f(\boldsymbol{W})
=
\frac{1}{N}(\boldsymbol{W}\boldsymbol{X}-\boldsymbol{Y})\boldsymbol{X}^\top+\lambda \boldsymbol{W}.
\]
In the derivative-free setting considered here, this gradient is assumed to be
unavailable; only function values \(f(\boldsymbol{W})\) can be evaluated.

The proposed derivative-free Muon framework can be applied directly. Given the
rank-one perturbation
\[
\boldsymbol{H}=\boldsymbol{u}\boldsymbol{v}^\top,
\qquad
\boldsymbol{u}\in\mathbb{R}^m,\quad \boldsymbol{v}\in\mathbb{R}^n,
\]
we evaluate
\[
f(\boldsymbol{W}+h \boldsymbol{u}\boldsymbol{v}^\top)
\]
and form the finite-difference score
\[
s(\boldsymbol{u},\boldsymbol{v})
=
\frac{f(\boldsymbol{W}+h \boldsymbol{u}\boldsymbol{v}^\top)-f(\boldsymbol{W})}{h}.
\]
For small \(h\), this quantity satisfies
\[
s(\boldsymbol{u},\boldsymbol{v})
\approx
\langle \nabla f(\boldsymbol{W}),\boldsymbol{u}\boldsymbol{v}^\top\rangle_F
=
\boldsymbol{u}^\top \nabla f(\boldsymbol{W})\boldsymbol{v}.
\]
Thus, rank-one perturbations probe bilinear information in the gradient
matrix without requiring explicit backpropagation.

The measured bilinear quantities provide the coefficients for the rank-one
surrogate used in Algorithm~\ref{alg:df_muon_lowrank}.

\subsection{Experimental Setup for the Single-Layer Example}

We consider the synthetic matrix regression model
\[
\boldsymbol{Y}=\boldsymbol{W}_\star \boldsymbol{X}+\boldsymbol{\varepsilon},
\]
where \(\boldsymbol{W}_\star\in\mathbb{R}^{m\times n}\) is the ground-truth matrix,
\(\boldsymbol{X}\in\mathbb{R}^{n\times N}\) is the data matrix, and \(\boldsymbol{\varepsilon}\) is a
noise matrix. The goal is to recover \(\boldsymbol{W}_\star\) by minimizing
\[
f(\boldsymbol{W})
=
\frac{1}{2N}\|\boldsymbol{W}\boldsymbol{X}-\boldsymbol{Y}\|_F^2
+
\frac{\lambda}{2}\|\boldsymbol{W}\|_F^2 .
\]

We compare the following methods:
\begin{itemize}
    \item Full-gradient Muon;
    \item entrywise derivative-free Muon;
    \item low-rank structured derivative-free Muon;
    \item basis-aligned rank-one derivative-free Muon (V2);
    \item direct structured search.
\end{itemize}
The performance is measured by the relative objective value
\[
\frac{f(\boldsymbol{W}_t)-f^\star}{f(\boldsymbol{W}_0)-f^\star},
\]
the relative recovery error
\[
\frac{\|\boldsymbol{W}_t-\boldsymbol{W}_\star\|_F}{\|\boldsymbol{W}_\star\|_F},
\]
and the number of function evaluations.

This experiment is well suited to comparing entrywise recovery with structure
recovery. Entrywise derivative-free Muon requires many function evaluations to
approximate the full \(m\times n\) gradient matrix, whereas the low-rank
structured method uses only rank-one probes and can be much cheaper when the
dominant gradient information is approximately low-rank.

\subsection{Experimental Protocol and Cases}

Every matrix-regression method is run for 200
iterations; the Muon variants use momentum parameter \(\beta=0.9\), forward-difference radius
\(h=10^{-4}\), and an SVD-based polar factor as the reference for
Muon orthogonalization.  The ridge optimum \(f^\star\) is computed analytically.
We use an economy-size SVD and return $UV^\top$ for
nonzero inputs; at rank deficiency this includes an orthogonal completion in
null directions. This convention differs from the compact nonzero-subspace
factor used in the sparse illustration above. Inputs with Frobenius norm at
most machine epsilon are mapped to zero. Randomness is fixed separately for each case.  Table~\ref{tab:case-settings}
summarizes the dimensions, target structure, noise levels, and algorithmic
parameters for the four cases.
The table lists the step sizes for the Muon variants. Direct search uses
$\eta=0.003$ in Cases 1--3 and $\eta=0.002$ in Case 4, with no momentum.
The case seeds are 1101, 2202, 3303, and 4404. For random low-rank DFO,
the surrogate is rescaled by $mn/r$; this fixed positive scaling
preserves the ideal polar direction.

\begin{table}[htbp]
\centering
\caption{Matrix-regression settings.  The number in the probe column
is the number of random rank-one probes used per iteration; \(q\) is the number
of direct-search candidates.}
\label{tab:case-settings}
\small
\begin{tabular}{cccccccccc}
\toprule
Case & \(m\) & \(n\) & \(N\) & target & noise & \(\lambda\) & \(\eta\) & probes & \(q\) \\
\midrule
1 & 8 & 8 & 256 & dense & 0.01 & \(10^{-4}\) & 0.035 & 4 & 8 \\
2 & 8 & 8 & 256 & dense & 0.10 & \(10^{-4}\) & 0.035 & 4 & 8 \\
3 & 8 & 8 & 256 & rank 2 & 0.01 & \(10^{-4}\) & 0.035 & 4 & 8 \\
4 & 16 & 16 & 512 & dense & 0.01 & \(10^{-4}\) & 0.025 & 8 & 12 \\
\bottomrule
\end{tabular}
\end{table}

We count every black-box objective evaluation.  With 200 updates,
full-gradient Muon records 201 objective values.  Entrywise DFO-Muon and V2 use
a cached base value, \(mn\) perturbations, and one evaluation of the
updated iterate per update, giving
\(1+200(1+mn)\) calls: 13,001 in Cases 1--3 and 51,401 in Case 4.  Random
low-rank DFO-Muon analogously costs \(1+200(1+r)\) calls, and direct search
costs \(1+200(1+q)\).  V2 therefore has the same query complexity as entrywise recovery.

\begin{figure}[htbp]
    \centering
    \includegraphics[width=0.49\textwidth]{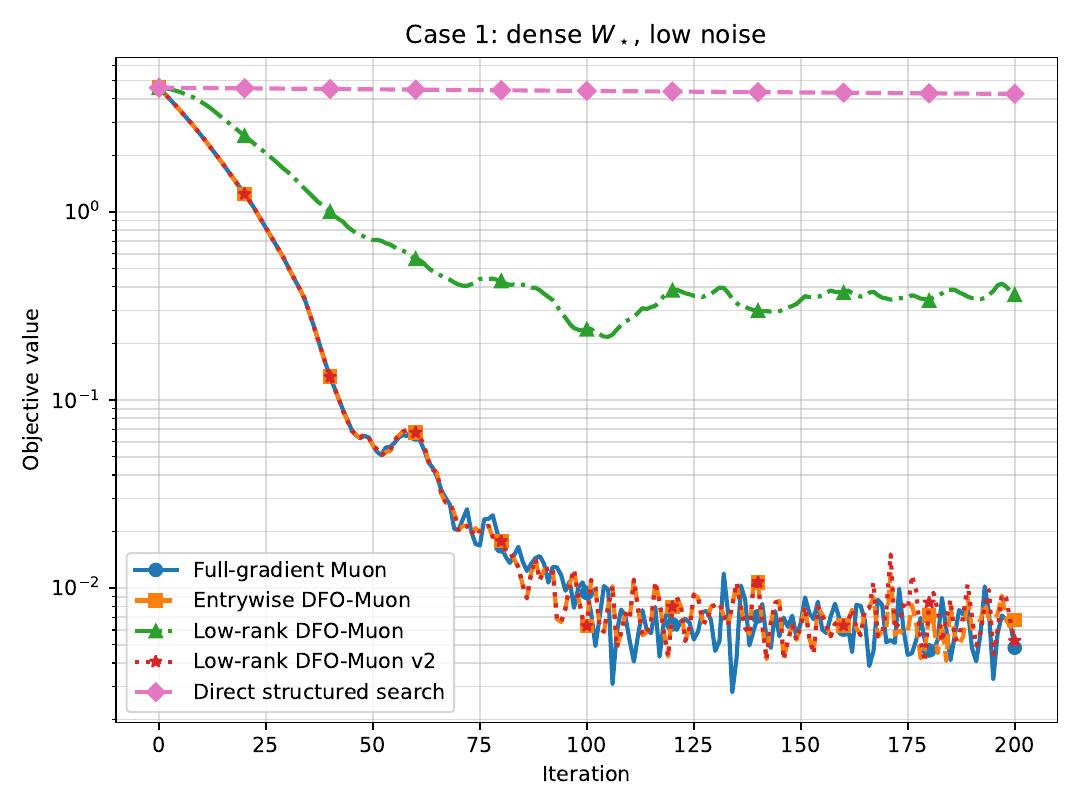}\hfill
    \includegraphics[width=0.49\textwidth]{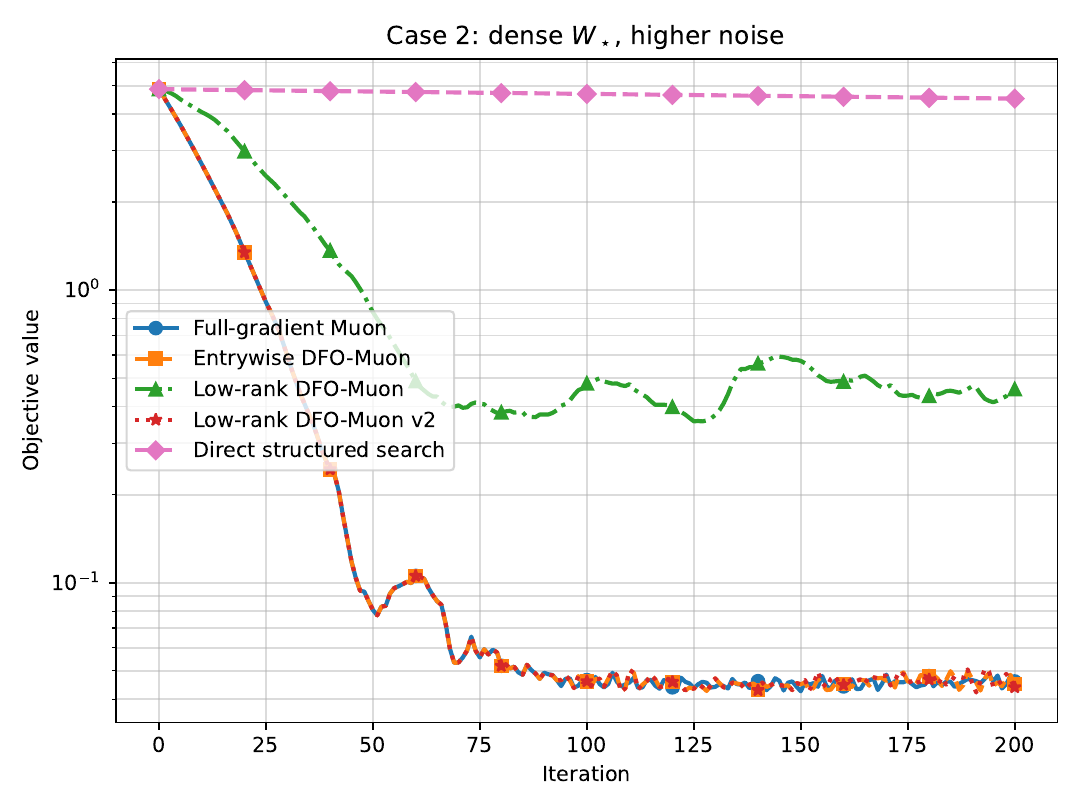}\\[0.5em]
    \includegraphics[width=0.49\textwidth]{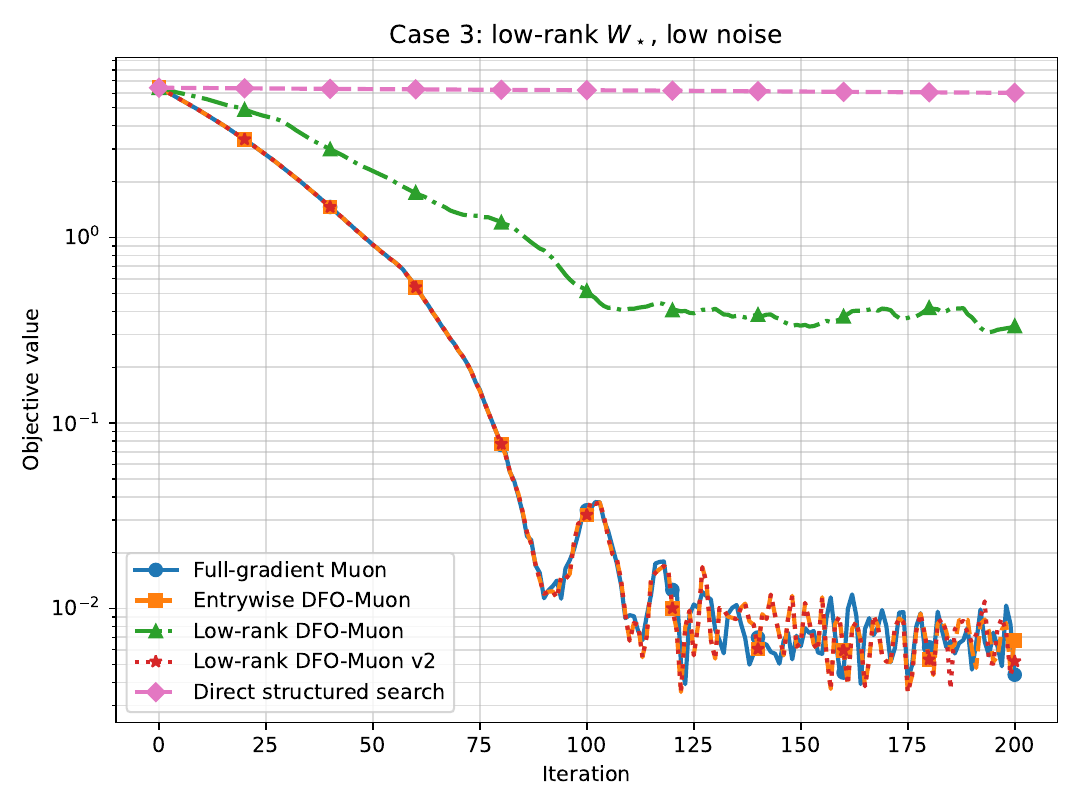}\hfill
    \includegraphics[width=0.49\textwidth]{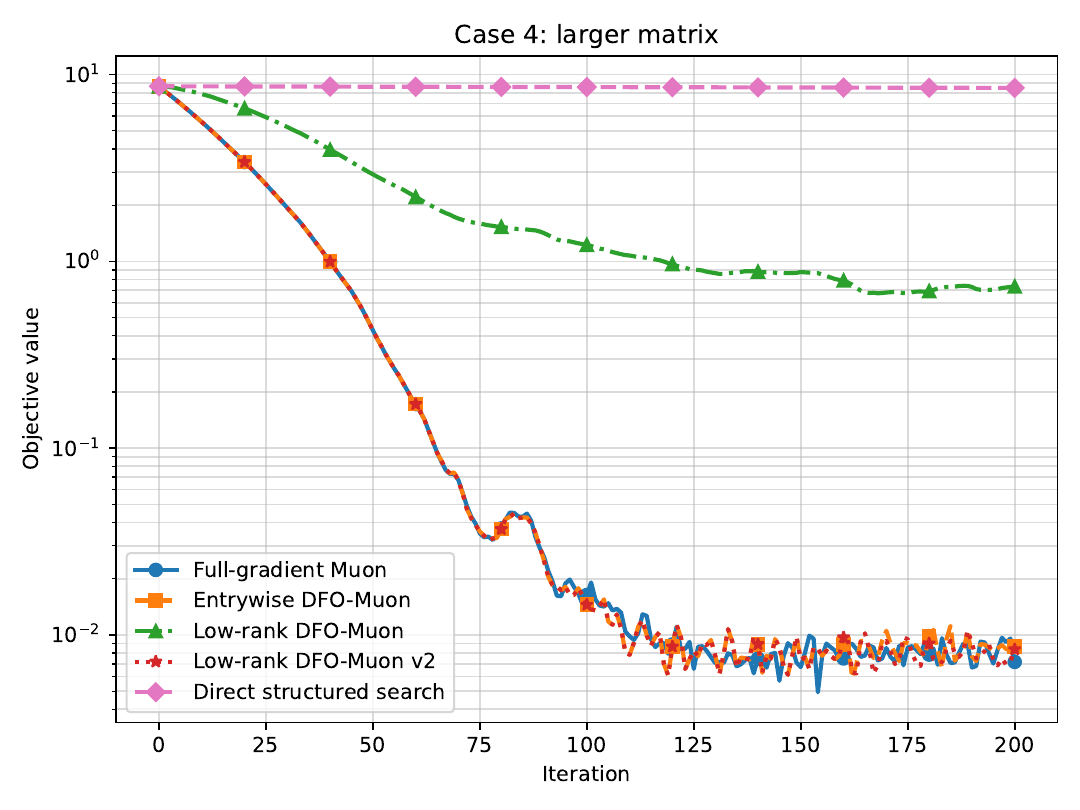}
    \caption{Objective value versus iteration for all four matrix-regression cases and
    all five methods.  Cases 1 and 2 use dense \(8\times8\) targets with low and
    higher noise, respectively; Case 3 uses a rank-2 \(8\times8\) target; Case 4
    uses a dense \(16\times16\) target.  The logarithmic vertical scale exposes
    both rapid convergence and terminal plateaus.}
    \label{fig:matrix-regression-all}
\end{figure}

\subsection{Terminal Results}

Table~\ref{tab:terminal-results} reports the final objective, normalized
objective gap, relative recovery error, and cumulative function evaluations.  Initial objectives are 4.578157,
4.864675, 6.419220, and 8.657680 for Cases 1--4, while the corresponding ridge
optima are 0.000835, 0.039036, 0.000997, and 0.001644.  The normalized gap is
therefore comparable across cases even though their raw objective scales and
irreducible noise floors differ.

\begin{table}[htbp]
\centering
\caption{Terminal results after 200 iterations.  ``Full'', ``Entry'',
``LR'', ``V2'', and ``Direct'' denote full-gradient Muon, entrywise DFO-Muon,
random low-rank DFO-Muon, basis-aligned low-rank DFO-Muon V2, and direct
structured search, respectively.}
\label{tab:terminal-results}
\scriptsize
\resizebox{\textwidth}{!}{%
\begin{tabular}{c l r r r r}
\toprule
Case & Method & Final \(f\) & Relative gap & Recovery error & Evaluations \\
\midrule
1 & Full   & 0.004820 & 0.000871 & 0.030115 & 201 \\
1 & Entry  & 0.006766 & 0.001296 & 0.035121 & 13,001 \\
1 & LR     & 0.362691 & 0.079054 & 0.286654 & 1,001 \\
1 & V2     & 0.005250 & 0.000964 & 0.032060 & 13,001 \\
1 & Direct & 4.245535 & 0.927333 & 0.964388 & 1,801 \\
\midrule
2 & Full   & 0.045994 & 0.001442 & 0.048530 & 201 \\
2 & Entry  & 0.045046 & 0.001245 & 0.041448 & 13,001 \\
2 & LR     & 0.459175 & 0.087064 & 0.313414 & 1,001 \\
2 & V2     & 0.043920 & 0.001012 & 0.038797 & 13,001 \\
2 & Direct & 4.516139 & 0.927774 & 0.964882 & 1,801 \\
\midrule
3 & Full   & 0.004402 & 0.000530 & 0.023890 & 201 \\
3 & Entry  & 0.006710 & 0.000890 & 0.031924 & 13,001 \\
3 & LR     & 0.333697 & 0.051837 & 0.236829 & 1,001 \\
3 & V2     & 0.005190 & 0.000653 & 0.026431 & 13,001 \\
3 & Direct & 6.026413 & 0.938798 & 0.970118 & 1,801 \\
\midrule
4 & Full   & 0.007174 & 0.000639 & 0.025408 & 201 \\
4 & Entry  & 0.008698 & 0.000815 & 0.028365 & 51,401 \\
4 & LR     & 0.734327 & 0.084644 & 0.294062 & 1,801 \\
4 & V2     & 0.008368 & 0.000777 & 0.027836 & 51,401 \\
4 & Direct & 8.486894 & 0.980270 & 0.990448 & 2,601 \\
\bottomrule
\end{tabular}}
\end{table}

\subsection{Interpretation of the Four Cases}

In Case 1, full-gradient Muon, entrywise DFO-Muon, and V2 rapidly reduce the
objective by roughly three orders of magnitude and finish with recovery errors
between 3.0\% and 3.5\%.  Random low-rank probing still produces substantial
descent: its objective falls from 4.578 to 0.363 using only 1,001 evaluations,
but its 28.7\% recovery error shows that four fresh random probes per iteration
do not identify enough of a dense gradient.  Direct structured search retains 92.7\% of the initial normalized gap.

Case 2 uses tenfold higher observation noise.  Its ridge
optimum rises to 0.039036, so the final raw objectives of the three exhaustive
or exact methods are necessarily higher.  V2 has the lowest recorded terminal
objective (0.043920) and recovery error (3.88\%), followed closely by entrywise
DFO-Muon and full-gradient Muon.  These small differences should not be read as
a systematic advantage over the analytic gradient: the oscillatory terminal
phase and a single fixed seed make endpoint rankings sensitive to the final
iterate.  Random low-rank probing plateaus at 0.459175, consistent with noisy,
sparsely sampled bilinear measurements.

Case 3 is consistent with a benefit from low-rank problem structure.  With a rank-2 target,
the random low-rank method obtains its best normalized gap (0.051837) and best
recovery error (23.7\%) of the four cases while retaining the 1,001-evaluation
budget.  Its error nevertheless remains larger than those of full-gradient Muon,
entrywise DFO-Muon, and V2. A single case does not establish how much of the
difference is attributable to target rank.
V2 again closely follows full-gradient Muon because its complete canonical basis
recovers essentially all coordinate information.

Case 4 increases the matrix dimension.  Increasing the parameter from 64 to 256
entries multiplies the query count of entrywise DFO-Muon and V2 from 13,001 to
51,401.  Their recovery errors remain small (2.84\% and 2.78\%), at the cost of exhaustive probing.  Random low-rank DFO-Muon uses
only 1,801 calls and reduces the objective from 8.658 to 0.734, a useful but
incomplete decrease.  The one-sided always-move direct search is least effective:
after 2,601 calls it retains 98.0\% of the initial normalized gap.

Across all cases, exhaustive finite differences and V2 reproduce high-quality
Muon directions, whereas random low-rank probes trade accuracy for 13--29 times
fewer evaluations than the exhaustive methods.  Direct random rank-one search
does not provide a competitive update under the present fixed step and
always-move rule.

\subsection{When Function Values Are More Reliable Than Gradients}
\label{subsec:noisy-gradient}

We next hold function evaluations exact and add noise to the gradient supplied
to Muon. This separates gradient-estimation error from noise in the objective.

\begin{figure}[htbp]
  \centering
  \includegraphics[width=\textwidth]{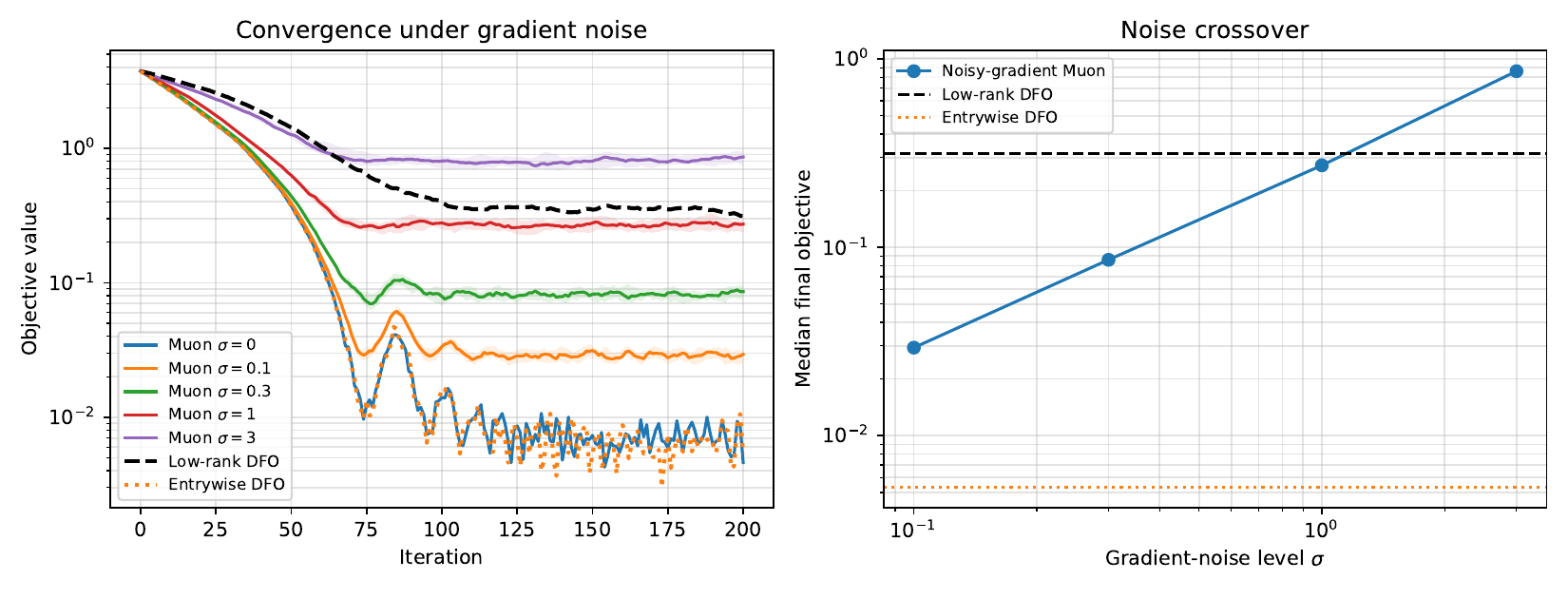}
  \caption{Comparison under an exact function-value oracle and an increasingly
  noisy gradient oracle.  Left: median objective trajectories, with interquartile
  bands for randomized runs.  Right: terminal-objective crossover.  Entrywise
  DFO is already more accurate at \(\sigma=0.1\); the cheaper four-probe low-rank
  DFO crosses noisy-gradient Muon between the tested levels \(\sigma=1\) and
  \(\sigma=3\).}
  \label{fig:noisy-gradient}
\end{figure}

Case 5 is the rank-2
\(8\times8\) matrix-regression problem with \(N=256\), observation-noise standard
deviation 0.01, and the same \(\lambda=10^{-4}\), \(\eta=0.035\), \(\beta=0.9\),
\(h=10^{-4}\), and 200 iterations as above.  Function evaluations are exact,
but the first-order oracle supplied to Muon is
\begin{equation}
\widetilde{\boldsymbol{G}}_t
=\nabla f(\boldsymbol{W}_t)+\boldsymbol{E}_t,
\qquad
(\boldsymbol{E}_t)_{ij}\overset{\mathrm{iid}}{\sim}
\mathcal{N}\!\left(0,
\sigma^2\frac{\|\nabla f(\boldsymbol{W}_0)\|_F^2}{mn}\right).
\label{eq:noisy-gradient-oracle}
\end{equation}
The noise variance is fixed at its initial scale. As the gradient shrinks,
the relative error therefore increases. This represents an estimator with a
persistent noise floor.  We test
\(\sigma\in\{0,0.1,0.3,1,3\}\) over 30 independent oracle-noise realizations.
Entrywise DFO is deterministic here, while random low-rank DFO is also repeated
30 times.  Shaded bands in Figure~\ref{fig:noisy-gradient} are interquartile
ranges, not confidence intervals.

\begin{table}[htbp]
\centering
\caption{Case 5 terminal results.  Quartiles summarize 30 runs except for the
deterministic entrywise method.  Function evaluations count objective queries;
the application-dependent cost of producing a noisy gradient is excluded.}
\label{tab:noisy-gradient}
\small
\begin{tabular}{lrrrrr}
\toprule
Method & \(\sigma\) & Median \(f_{200}\) & Q1 & Q3 & Evaluations \\
\midrule
Noisy-gradient Muon & 0   & 0.004601 & 0.004601 & 0.004601 & 201 \\
Noisy-gradient Muon & 0.1 & 0.029352 & 0.025675 & 0.031733 & 201 \\
Noisy-gradient Muon & 0.3 & 0.085930 & 0.075410 & 0.096034 & 201 \\
Noisy-gradient Muon & 1   & 0.272628 & 0.237528 & 0.290655 & 201 \\
Noisy-gradient Muon & 3   & 0.862138 & 0.740106 & 0.954609 & 201 \\
Low-rank DFO-Muon   & --- & 0.313504 & 0.265583 & 0.350586 & 1,001 \\
Entrywise DFO-Muon  & --- & 0.005336 & 0.005336 & 0.005336 & 13,001 \\
\bottomrule
\end{tabular}
\end{table}

With an exact gradient, Muon has the lowest terminal objective: its median final objective is 0.004601, slightly
below entrywise DFO. Forward-difference truncation error is one possible
contributor to this difference.  Once
\(\sigma=0.1\), however, entrywise DFO reaches 0.005336 while noisy-gradient
Muon stops at 0.029352.  Thus a sufficiently accurate function-value oracle can
outperform a mildly noisy first-order oracle, although exhaustive coordinate
probing costs roughly 65 times as many objective calls per run.  At \(\sigma=3\) low-rank DFO reaches a median
0.313504, substantially below Muon's 0.862138, with 1,001 rather than 13,001
function evaluations.  The measured crossover for low-rank DFO lies between
the sampled levels 1 and 3; the experiment does not claim a universal threshold.

The comparison measures objective accuracy and query counts, not elapsed
time. A runtime comparison would require the cost of the gradient oracle,
including any additional samples needed to reduce its variance. The present
experiment uses synthetic gradient corruption and does not measure these costs.

\subsection{Neural-Network Matrix Optimization with a Noisy Backward Pass}
\label{subsec:neural-noisy-gradient}

The next experiment trains the hidden-layer matrix of a teacher--student
neural network,
\begin{equation}
 \widehat{\boldsymbol{Y}}(\boldsymbol{W})
 =\boldsymbol{A}\tanh(\boldsymbol{W}\boldsymbol{X}),
 \qquad
 \boldsymbol{Y}=\boldsymbol{A}\tanh(\boldsymbol{W}_\star\boldsymbol{X}),
 \label{eq:fixed-readout-network}
\end{equation}
where the readout matrix \(\boldsymbol{A}\in\mathbb{R}^{4\times8}\) is fixed and
only the hidden-layer weight
\(\boldsymbol{W}\in\mathbb{R}^{8\times6}\) is trained.  For \(N=512\) input
samples, the loss is
\[
 f(\boldsymbol{W})=
 \frac{1}{2N}\left\|
 \boldsymbol{A}\tanh(\boldsymbol{W}\boldsymbol{X})-\boldsymbol{Y}
 \right\|_F^2.
\]
The forward computation requires only two matrix multiplications and an
elementwise \(\tanh\).  Its exact backpropagated gradient is
\begin{equation}
 \nabla f(\boldsymbol{W})=
 \frac{1}{N}\left[
 \boldsymbol{A}^{\top}\boldsymbol{R}\odot
 \left(\boldsymbol{1}-\boldsymbol{H}\odot\boldsymbol{H}\right)
 \right]\boldsymbol{X}^{\top},
 \quad
 \boldsymbol{H}=\tanh(\boldsymbol{W}\boldsymbol{X}),\quad
 \boldsymbol{R}=\boldsymbol{A}\boldsymbol{H}-\boldsymbol{Y}.
 \label{eq:network-gradient}
\end{equation}
This exact expression is used only to create a controlled reference and a noisy
backward oracle.  In a black-box implementation, the same forward map could be
an external simulator, quantized accelerator, discrete activation pipeline, or
hardware measurement for which Equation~\eqref{eq:network-gradient} is not
available to the optimizer.

At every iteration, standard Muon forms momentum from either the exact gradient
or the corrupted gradient in Equation~\eqref{eq:noisy-gradient-oracle}, and
updates the hidden weight matrix with its polar factor.  The two derivative-free
variants instead perturb this same \(8\times6\) hidden matrix and evaluate the
network loss: entrywise DFO uses all 48 canonical rank-one directions, while
low-rank DFO uses six random normalized outer products per iteration.  All
methods start from \(\boldsymbol{W}_0=0\) and use 300 iterations,
\(\eta=0.025\), \(\beta=0.9\), and \(h=10^{-4}\).  Random methods are repeated
over 30 seeds.

\begin{figure}[htbp]
 \centering
 \includegraphics[width=\textwidth]{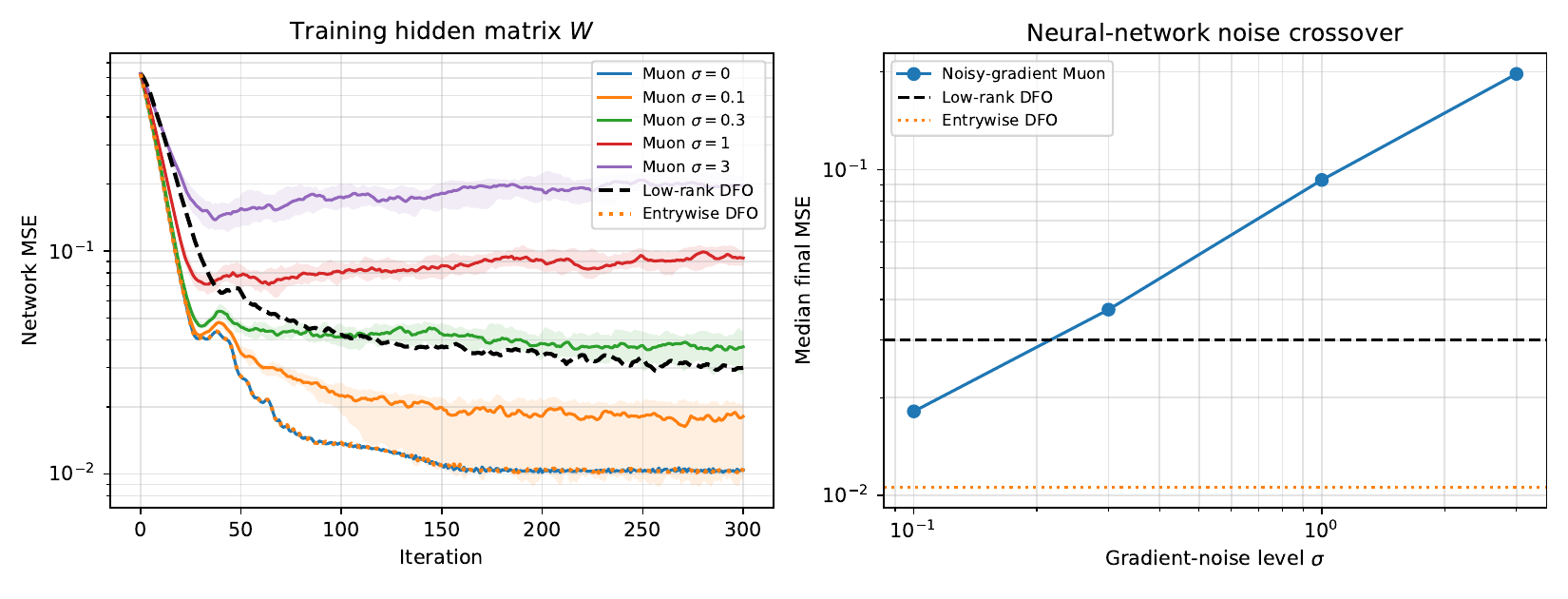}
 \caption{Training the hidden matrix of the fixed-readout neural network in
 Equation~\eqref{eq:fixed-readout-network}.  Left: median training-loss
 trajectories and interquartile bands.  Right: terminal loss against the noise
 level of the backward gradient.  Low-rank DFO-Muon overtakes noisy-gradient
 Muon between the tested levels \(\sigma=0.1\) and \(\sigma=0.3\).}
 \label{fig:neural-noisy-gradient}
\end{figure}

\begin{table}[htbp]
\centering
\caption{Terminal neural-network loss after 300 hidden-matrix updates.  The
quartiles summarize 30 runs except for deterministic entrywise DFO-Muon.}
\label{tab:neural-noisy-gradient}
\small
\begin{tabular}{lrrrrr}
\toprule
Method & \(\sigma\) & Median \(f_{300}\) & Q1 & Q3 & Evaluations \\
\midrule
Noisy-gradient Muon & 0   & 0.010430 & 0.010430 & 0.010430 & 301 \\
Noisy-gradient Muon & 0.1 & 0.018142 & 0.009370 & 0.020764 & 301 \\
Noisy-gradient Muon & 0.3 & 0.037245 & 0.028105 & 0.044499 & 301 \\
Noisy-gradient Muon & 1   & 0.093108 & 0.086176 & 0.099475 & 301 \\
Noisy-gradient Muon & 3   & 0.196978 & 0.172462 & 0.217890 & 301 \\
Low-rank DFO-Muon   & --- & 0.030017 & 0.025091 & 0.036054 & 2,101 \\
Entrywise DFO-Muon  & --- & 0.010614 & 0.010614 & 0.010614 & 14,701 \\
\bottomrule
\end{tabular}
\end{table}

With an exact backward pass, Muon attains the lowest loss, 0.010430, and
entrywise DFO nearly reproduces it at 0.010614.  At \(\sigma=0.1\), entrywise
DFO already improves on the median noisy-gradient Muon endpoint, although the
overlapping run distribution cautions against interpreting this as a sharp
threshold.  At \(\sigma=0.3\), six-probe low-rank DFO reaches 0.030017 and
overtakes noisy-gradient Muon at 0.037245.  For \(\sigma=1\) and 3, the gap
widens substantially: noisy-gradient Muon ends at 0.093108 and 0.196978,
respectively, while low-rank DFO is unaffected because its oracle consists of
exact forward losses.

Low-rank DFO uses 2,101 forward calls, compared with 301 recorded losses for
Muon and 14,701 calls for entrywise DFO. Its lower loss under sufficiently
large gradient noise comes with this additional query cost. Whether it saves
computation depends on the cost of obtaining a more accurate backward pass.

\subsection{Black-Box Policy Optimization on CartPole}
\label{subsec:cartpole-dfo}

We finally test the proposed optimizer in a setting where the objective is an
episode score and no derivative through the environment is supplied.  A
stochastic policy with two matrix-valued layers is used,
\begin{equation}
 \boldsymbol{z}=\tanh(\boldsymbol{W}_1\boldsymbol{s}),\qquad
 \pi_{\boldsymbol{W}}(a\mid\boldsymbol{s})
 =\operatorname{softmax}(\boldsymbol{W}_2\boldsymbol{z})_a,
 \quad
 \boldsymbol{W}_1\in\mathbb{R}^{8\times4},\quad
 \boldsymbol{W}_2\in\mathbb{R}^{2\times8}.
 \label{eq:cartpole-policy}
\end{equation}
The black-box loss is the negative mean CartPole return,
\(f(\boldsymbol{W})=-\mathbb{E}[R(\boldsymbol{W})]\), with a horizon of 500
steps.  The optimizer receives only complete episode returns; the derivative-free methods differentiate neither the simulator
dynamics nor the action probabilities.

At update \(t\), each probe independently constructs an isotropic rank-one
direction for each layer,
\[
 \boldsymbol{H}_{\ell,k}
 =\sqrt{m_\ell n_\ell}\,
 \frac{\boldsymbol{u}_{\ell,k}}{\|\boldsymbol{u}_{\ell,k}\|_2}
 \frac{\boldsymbol{v}_{\ell,k}^{\top}}{\|\boldsymbol{v}_{\ell,k}\|_2},
 \qquad \ell\in\{1,2\}.
\]
Both layers are perturbed simultaneously and the loss directional derivative is
estimated by
\begin{equation}
 y_{t,k}=\frac{\widehat f(\boldsymbol{W}_t+h\boldsymbol{H}_k)
 -\widehat f(\boldsymbol{W}_t-h\boldsymbol{H}_k)}{2h},
 \qquad
 \widehat{\boldsymbol{G}}_{\ell,t}
 =\frac{1}{q}\sum_{k=1}^{q}y_{t,k}\boldsymbol{H}_{\ell,k}.
 \label{eq:cartpole-zo-estimator}
\end{equation}

In addition to independent probes, we test an enhanced construction that makes
the probes within each update mutually orthogonal in the Frobenius inner
product.  For each layer, independent Gaussian matrices are orthogonalized to
form bases $\boldsymbol{Q}_L=[\boldsymbol{q}_{L,1},\ldots]$ and
$\boldsymbol{Q}_R=[\boldsymbol{q}_{R,1},\ldots]$.  Distinct index pairs are
sampled without replacement and
\begin{equation}
 \boldsymbol{H}_k=\sqrt{mn}\,
 \boldsymbol{q}_{L,i_k}\boldsymbol{q}_{R,j_k}^{\top},
 \qquad
 \langle\boldsymbol{H}_k,\boldsymbol{H}_{k'}\rangle_F=0
 \quad(k\ne k').
 \label{eq:orthogonal-rank-one-probes}
\end{equation}
The random rotations preserve isotropic marginal directions, while sampling
without replacement ensures orthogonality of the selected matrix measurements.  We call the resulting method enhanced DFO-Muon.  A preliminary
ablation also tested truncating the polar factor to the singular subspace that
contained 90\% of the momentum energy.  That modification slowed learning in
this control problem, so the reported enhanced method retains the full polar
factor in the main comparison; the ablation is exploratory.
For every positive--negative pair, the same initial states and the same uniform
random numbers for action sampling are reused.  This common-random-number
coupling reduces variance without exposing a pathwise derivative.  DFO-Muon
forms layerwise momentum from Equation~\eqref{eq:cartpole-zo-estimator} and
applies a polar factor to each layer.  The query-matched baseline uses the same
rank-one probes and momentum but applies a single normalized zeroth-order SGD
step.  Thus the comparison changes the matrix update geometry, not the
black-box observations.

We also include a standard Muon baseline driven by a likelihood-ratio policy
gradient.  It uses the identity
\begin{equation}
 \nabla_{\boldsymbol{W}}J(\boldsymbol{W})
 =\mathbb{E}\left[\sum_{t=0}^{T-1}
 G_t\nabla_{\boldsymbol{W}}\log
 \pi_{\boldsymbol{W}}(a_t\mid\boldsymbol{s}_t)\right],
 \label{eq:reinforce-muon}
\end{equation}
where $J(\boldsymbol{W})=\mathbb{E}[R(\boldsymbol{W})]$ is the expected
return and $G_t$ is the reward-to-go.  The policy-gradient estimator uses a time-indexed batch
baseline and normalizes the sampled advantages before applying layerwise Muon
momentum and polar factors.  Equation~\eqref{eq:reinforce-muon} differentiates
the policy log probability, but does not differentiate the CartPole transition
map or reward.  We refer to this method as policy-gradient Muon.

The original zeroth-order methods use \(q=8\) antithetic probes, three training
episodes per perturbed policy, and \(h=0.1\).  Enhanced DFO-Muon reallocates the
same budget to 12 orthogonal probes, two episodes per perturbed policy, and
\(h=0.2\).  Thus every zeroth-order update consumes 48 episodes.  Policy-gradient
Muon uses a batch of 48 episodes per update as well, so all four methods have
the same episode budget at every plotted iteration.  All methods use
\(\eta=0.1\), \(\beta=0.8\), and 70 policy updates.
All results are averaged over five fixed random seeds.  A separate fixed set of
40 episodes per seed is used only to report each policy's mean return and is not
included in the optimizer query count.  A run is called solved at its first
reported mean return of at least 475.

\begin{figure}[htbp]
 \centering
 \includegraphics[width=\textwidth]{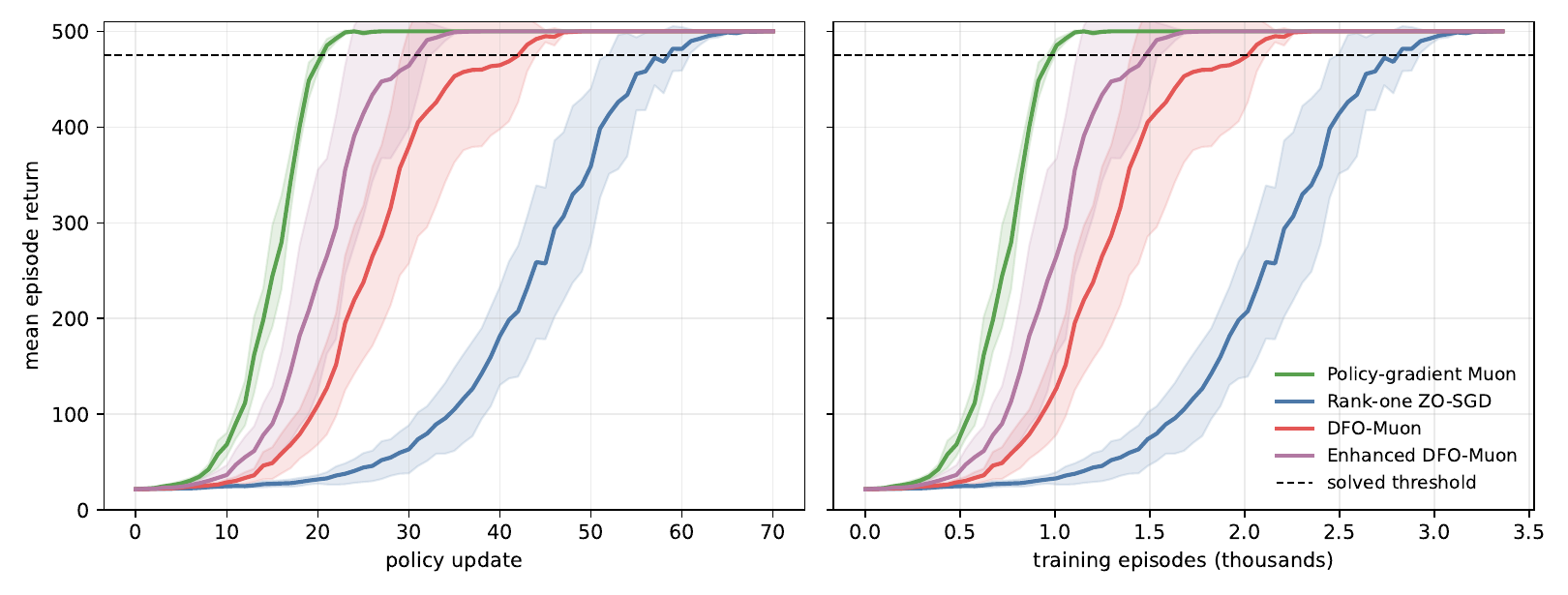}
 \caption{CartPole policy optimization.  Curves show the five-seed mean and
 shaded bands show one standard deviation, against updates (left) and training
 episodes (right).  The original zeroth-order pair uses identical antithetic
 rank-one measurements; the enhanced variant replaces independent directions
 with orthogonal rank-one probes.  Policy-gradient Muon uses REINFORCE and no
 environment derivative.}
 \label{fig:cartpole-dfo}
\end{figure}

\begin{table}[htbp]
\centering
\caption{CartPole results over five seeds.  Every method uses 48 training
episodes per update.  The solved threshold is a reported mean return of 475.}
\label{tab:cartpole-dfo}
\small
\begin{tabular}{lrrrr}
\toprule
Method & Final return & Solved seeds & First solved update & Episodes to solve \\
\midrule
Policy-gradient Muon & $500.0\pm0.0$ & 5/5 & 20.4 & 979.2 \\
Rank-one ZO-SGD     & $500.0\pm0.0$ & 5/5 & 56.4 & 2,707.2 \\
DFO-Muon            & $500.0\pm0.0$ & 5/5 & 34.2 & 1,641.6 \\
Enhanced DFO-Muon   & $500.0\pm0.0$ & 5/5 & 27.2 & 1,305.6 \\
\bottomrule
\end{tabular}
\end{table}

All four methods reach the maximum return in all five runs.  Policy-gradient
Muon is strongest, crossing 475 after 20.4 updates and 979.2 episodes on average.
Enhanced DFO-Muon is the best black-box method, requiring 27.2 updates and
1,305.6 episodes, compared with 34.2 and 1,641.6 for the original DFO-Muon and
56.4 and 2,707.2 for rank-one ZO-SGD.  Mean return over the complete trajectory
is 395.8, 358.4, 320.6, and 203.1, respectively.  With shared rank-one measurements, DFO-Muon solves this task in fewer
training episodes than normalized ZO-SGD. The enhanced variant requires fewer
episodes still, although it changes the probe count, episodes per probe, and
perturbation radius along with probe orthogonality. Their individual effects
cannot be separated by this comparison. Policy-gradient Muon requires the
fewest episodes of the four methods.

These results are limited to a small policy and five seeds. Episode lengths
vary, so equal episode budgets need not imply equal numbers of environment
steps. The enhanced settings were selected through an exploratory ablation.

\subsubsection{Crossover under a noisy policy gradient}

The preceding comparison assumes that the REINFORCE gradient is available at
its native sampling accuracy.  To model an unreliable backward or distributed
gradient pipeline, we additionally corrupt each sampled layer gradient by
\begin{equation}
 \widetilde{\boldsymbol{G}}_{\ell,t}
 =\boldsymbol{G}_{\ell,t}
 +\sigma\frac{\|\boldsymbol{G}_{\ell,t}\|_F}{\sqrt{m_\ell n_\ell}}
 \boldsymbol{Z}_{\ell,t},
 \qquad
 (\boldsymbol{Z}_{\ell,t})_{ij}\overset{\mathrm{iid}}{\sim}\mathcal{N}(0,1).
 \label{eq:noisy-policy-gradient}
\end{equation}
This is relative rather than fixed-floor corruption: $\sigma=1$ makes the
entrywise noise standard deviation equal to the RMS of the sampled gradient in
that layer.  Noisy policy-gradient Muon otherwise uses the same REINFORCE
estimator, momentum, learning rate, initialization, and 48-episode update
budget.  We test $\sigma\in\{0,1,3,10\}$ over the same five seeds.

\begin{figure}[htbp]
 \centering
 \includegraphics[width=\textwidth]{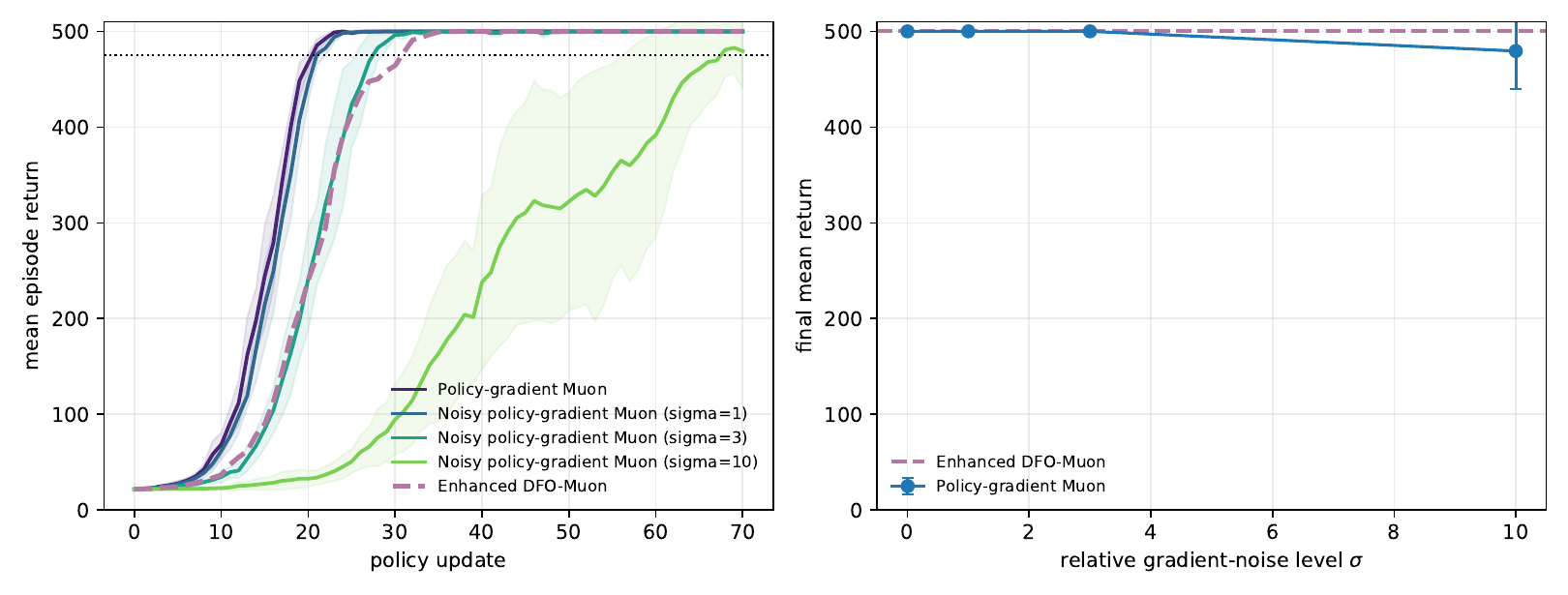}
 \caption{Policy-gradient Muon under the controlled corruption in
 Equation~\eqref{eq:noisy-policy-gradient}.  Left: mean return trajectories with
 one-standard-deviation bands and the enhanced DFO-Muon reference.  Right:
 final return against relative gradient-noise level; error bars show one
 standard deviation.}
 \label{fig:cartpole-noisy-policy-gradient}
\end{figure}

\begin{table}[htbp]
\centering
\caption{CartPole crossover under noisy policy gradients.  The mean time to
solve is computed over successful seeds; every update uses 48 episodes.}
\label{tab:cartpole-noisy-policy-gradient}
\small
\begin{tabular}{lrrrr}
\toprule
Method & $\sigma$ & Final return & Solved seeds & Episodes to solve \\
\midrule
Policy-gradient Muon       & 0  & $500.0\pm0.0$  & 5/5 & 979.2 \\
Noisy policy-gradient Muon & 1  & $500.0\pm0.0$  & 5/5 & 1,036.8 \\
Noisy policy-gradient Muon & 3  & $500.0\pm0.0$  & 5/5 & 1,296.0 \\
Noisy policy-gradient Muon & 10 & $479.4\pm44.3$ & 4/5 & 2,796.0 \\
Enhanced DFO-Muon          & ---& $500.0\pm0.0$  & 5/5 & 1,305.6 \\
\bottomrule
\end{tabular}
\end{table}

At $\sigma=1$, policy-gradient Muon retains a clear advantage.  At
$\sigma=3$, it and enhanced DFO-Muon are effectively tied in this resolution,
requiring 1,296.0 and 1,305.6 episodes on average.  At $\sigma=10$, enhanced
DFO-Muon is more reliable in these five runs: all five runs solve and finish at 500,
whereas noisy policy-gradient Muon solves four runs and has mean final return
479.4 with standard deviation 44.3.  The measured crossover therefore lies
above the mild-noise regime and between the sampled levels 3 and 10; it is not
a universal threshold, because the corruption is synthetic and scales with the
current sampled-gradient RMS.

\subsection{Scope and Reproducibility Limitations}

The numerical results use the dimensions, seeds, step sizes, and probe budgets
specified above. Orthogonalization is performed using an SVD polar factor;
the results do not establish identical behavior for a finite Newton--Schulz
iteration. The matrix-regression cases use one fixed seed each, and the
CartPole study uses five seeds and a small policy. The available materials do
not suffice to independently regenerate the CartPole runs. These experiments
therefore provide limited evidence of performance beyond the tested settings.
The paper provides no general convergence analysis of the proposed methods.

\section{Conclusion}

Structured finite differences provide a way to construct Muon updates from
function values. Exhaustive coordinate probing closely follows full-gradient
Muon in the matrix-regression experiments, while random rank-one probing uses
13--29 times fewer evaluations and gives larger terminal errors. Exhaustive
basis-aligned V2 is a scaled form of coordinate finite differences and has the
same query cost. Direct rank-one search performs poorly with the fixed-step,
always-move rule tested here.

The noisy-gradient experiments show a tradeoff between gradient reliability
and the cost of additional function evaluations. In CartPole, enhanced
DFO-Muon uses fewer episodes than the other tested zeroth-order methods, but
more than policy-gradient Muon with uncorrupted gradients. It solves all five
runs at the highest tested corruption level, where noisy policy-gradient Muon
solves four. These comparisons use synthetic gradient noise and small problems;
they do not establish a runtime advantage in larger applications.

Further work is needed to relate probe selection to the quality of the polar
update and to obtain convergence guarantees. Adaptive query budgets and step
acceptance rules may help control the errors introduced by partial gradient
information.

\bibliographystyle{plainnat}
\bibliography{references}

\end{document}